\documentclass[11pt]{article}

\usepackage{amsfonts,amsthm,amsmath,latexsym,amssymb}
\usepackage[dvips]{graphics}
\usepackage[dvips]{graphics}

\input epsf.tex
\newdimen\epsfxsize
\newdimen\epsfysize

\renewcommand{\a}{\alpha}
\renewcommand{\t}{\tau}
\newcommand{\C}{\mathbb C}
\renewcommand{\H}{\mathbb H}

\newcommand{\R}{\mathbb R}

\newcommand{\diam}{\text{diam }}

\renewcommand{\Im}{\operatorname{Im}}
\renewcommand{\Re}{\operatorname{Re}}

\newtheorem{thm}{Theorem}[section]
\newtheorem{prop}[thm]{Proposition}

\newtheorem{corollary}[thm]{Corollary}

\newtheorem{lemma}[thm]{Lemma}

\makeatletter

\@addtoreset{equation}{section} \makeatother

\newcommand{\karo}{\hfill$\Box$}

\def\wt{\widetilde}
\def\wh{\widehat}
\def\E{{\mathbb E}}
\def\P{{\mathbb P}}
\def\LL{{\cal L}}
\def\eps{\varepsilon}
\def\pf{\noindent {\bf Proof}. }
\def\1{{\bf 1}}

\begin{document}

%\baselineskip=1.2\baselineskip
\begin{doublespace}

\title{\bf Schramm-Loewner Equations Driven by Symmetric Stable Processes}
\bigskip
\author{{\bf Zhen-Qing Chen}\footnote{Research supported in part by NSF
Grant DMS-0600206.}  \quad and \quad {\bf Steffen
Rohde}\footnote{Research supported in part by NSF Grants DMS-0501726
and DMS-0244408.}
\medskip \\
Department of Mathematics\\
University of Washington\\
Seattle, WA 98195, USA \medskip\\
\texttt{zchen@math.washington.edu} \ and \ {\texttt
{rohde@math.washington.edu }}}

\date{(August 9, 2007)}
\maketitle

\abstract{We consider shape, size and regularity of the hulls $K_t$
of the chordal Schramm-Loewner evolution driven by a symmetric
$\a$-stable process. We obtain derivative estimates, show that the
domains $\H\setminus K_t$ are H\"older domains, prove that $K_t$ has
 Hausdorff dimension 1,
and show that the trace is right-continuous with left limits
almost surely.}

\section{Introduction and Results}\label{s0}

The Loewner differential equation (LE for short)
\begin{equation}\label{eqn:1.1}
\partial_t g_t(z) = \frac2{g_t(z)-W_t} \ ,\quad g_0(z)=z
\end{equation}
takes as input a real-valued function $W_t$ ($t\geq0$) and produces
an increasing family of sets $(K_t)_{t\geq0}$ such that $g_t$ is the
(suitably normalized) conformal map from $\H\setminus K_t$ onto the
upper halfplane $\H.$ See Section \ref{s2}. The Schramm Loewner
Evolution $SLE_{\kappa}$ is the random process $K_t$ (or $g_t$) when
$W_t = B_{\kappa t}$ where $B_t$ is Brownian motion. See \cite{S}
and the references therein.

The spectacular success of $SLE_{\kappa}$
in describing scaling limits of
lattice models and in resolving numerous questions from probability
and mathematical physics motivates the study of the Loewner equation
driven by other stochastic processes. Roughly speaking, if the
driving function is sufficiently continuous, then LE produces a
continuous curve $\gamma(t)\in\overline{\H}$ defined by
$g_t(\gamma(t))=W_t.$ This so-called {\it trace} generates the hull
in the sense that $K_t=\gamma[0,t]$ (if $\gamma$ is not a simple
curve, one has to add the filled-in loops). If $W$ has a
discontinuity at time $t$, then $\gamma$ has a discontinuity too and
the trace grows a "branch". In fact, if $W$ is piecewise constant,
then $K$ is a union of analytic curves (and the $n$-th of these
curves is a geodesic for the hyperbolic metric in the half plane
minus the previous $n-1$ curves). Thus tree-like sets $K$ can be
described by LE with discontinuous driving term. In the mathematical
physics literature, the LE driven by the symmetric $\alpha$-stable
process $S_t$ (plus Brownian motion) has first appeared in
\cite{ROKG}. A mathematically rigorous treatment of some elementary
properties is in \cite{GW}.

Another motivation for studying random families of conformal maps
comes from a circle of problems known in the complex analysis
literature as Brennan's conjecture,  see \cite{B} or \cite{P}. The
problem is to maximize
$$
\beta_f(p) = \limsup_{r\to 1}\ \frac{\log \bigl(\int_0^{2\pi}
|f'(r e^{it})|^p dt\bigr)}{|\log(1-r)|}
$$
 over all bounded conformal
maps $f$ of the unit disc. While it is conjectured that
$\beta(p):=\sup_f \beta_f(p) = p^2/4$ for $-2\leq p\leq 2,$ there is
no proof of either $\beta(p)\leq p^2/4$ or $\beta(p)\geq p^2/4$, for
any nontrivial value of $p.$ The lower bound just requires one
example $f,$ but there are no candidates for extreme domains. From
work of Carleson, Jones, Makarov and others it is know that
extremals can be found amongst domains with self-similar boundary,
and that extremal boundaries can be approximated by "dendrites".
Whereas it is difficult to compute the above integral means for
individual functions $f$, it could be easier to estimate the
expected value
$$   \E \left[ \int_0^{2\pi}|f'(r e^{it})|^p dt \right]
$$
because in a rotationally invariant family this amounts to computing $\E[|f'(r)|^p]$.
The computations in \cite{RS} showed that
Brownian SLE
does not produce examples close to extremal.
At the 2001/02 Mittag-Leffler program "Probability and Conformal Mappings",
Nikolai Makarov
and the second author tried to find stochastic processes that produced large
integral
means, and recognized that it would be interesting to study LE driven by
symmetric stable processes.
The second author would like to thank Nick for these stimulating conversations. In 2003, Daniel Meyer
(then graduate student at University of Washington, Seattle) performed computer experiments that suggested
a nontrivial and perhaps even close to extremal integral means spectrum for the stable LE.

In this paper, we will consider LE driven by symmetric
$\alpha$-stable process $W_t=S_t,$ see Section \ref{s1} for the
definition of symmetric stable processes and some of the basic
properties. As $W_t$ satisfies a scaling relation different from
Brownian scaling, stable LE does not exhibit scale invariance and
thus it is no surprise that rescaling the hulls leads to
deterministic sets. Indeed,
we show in Section 3 that for $0<\a<2$,   as $s \to0,$ the rescaled
hulls $\frac1{s} K_{s^2}$ converge to the vertical line segment
$[0,2i]$ (in the Hausdorff metric) in probability. On the other
hand, for all $\eps>0,$
$$ \lim_{s\to \infty}
\P\left( \frac1s K_{s^2} \cap \{y>\eps\} \not= \emptyset \right) = 0
.
$$

\bigskip
\noindent We will then consider continuity  and metric properties of
the hulls by analyzing the backward flow
\begin{equation}\label{eqn:1.2}
\partial_t f_t(z) = - \frac2{f_t(z)-W_t} \ ,\quad f_0(z)=z.
\end{equation}
 For each fixed $t>0$, this random
conformal map $f_t(z)$ of $\H$ has the same distribution as
$g_t^{-1}(z-W_t)+W_t$ and thus $K_t$ has the same distribution as
$\H\setminus f_t(\H)-W_t.$
However as a family of maps, $\{f_t(\cdot ), t\geq 0\}$ does not have
the same distribution as $\{g^{-1}_t(\cdot ), t\geq 0\}$
(see the discussion at the beginning of Section \ref{Derivative Estimates}).
Write
$$ f_t(z)-W_t=X_t+iY_t, \qquad t\geq 0.
$$
It is easy to see that $Y_t$ is increasing in $t\geq 0$. We prove in
Section 4 that for $z=x+iy\in\H$ with $y<1$, if $\alpha\in[1,2)$,
$Y$ reaches height 1 almost surely when $\alpha \in [1, 2)$, and $Y$
does not reach height 1 with positive probability when $\alpha \in
(0, 1)$.

Below are some computer simulations for SLE driven by Cauchy stable
processes, with $t=0.1, 1, 10$ and $t=100$ respectively.

\bigskip \vbox{ \epsfxsize=3.0in
  \centerline{\epsffile{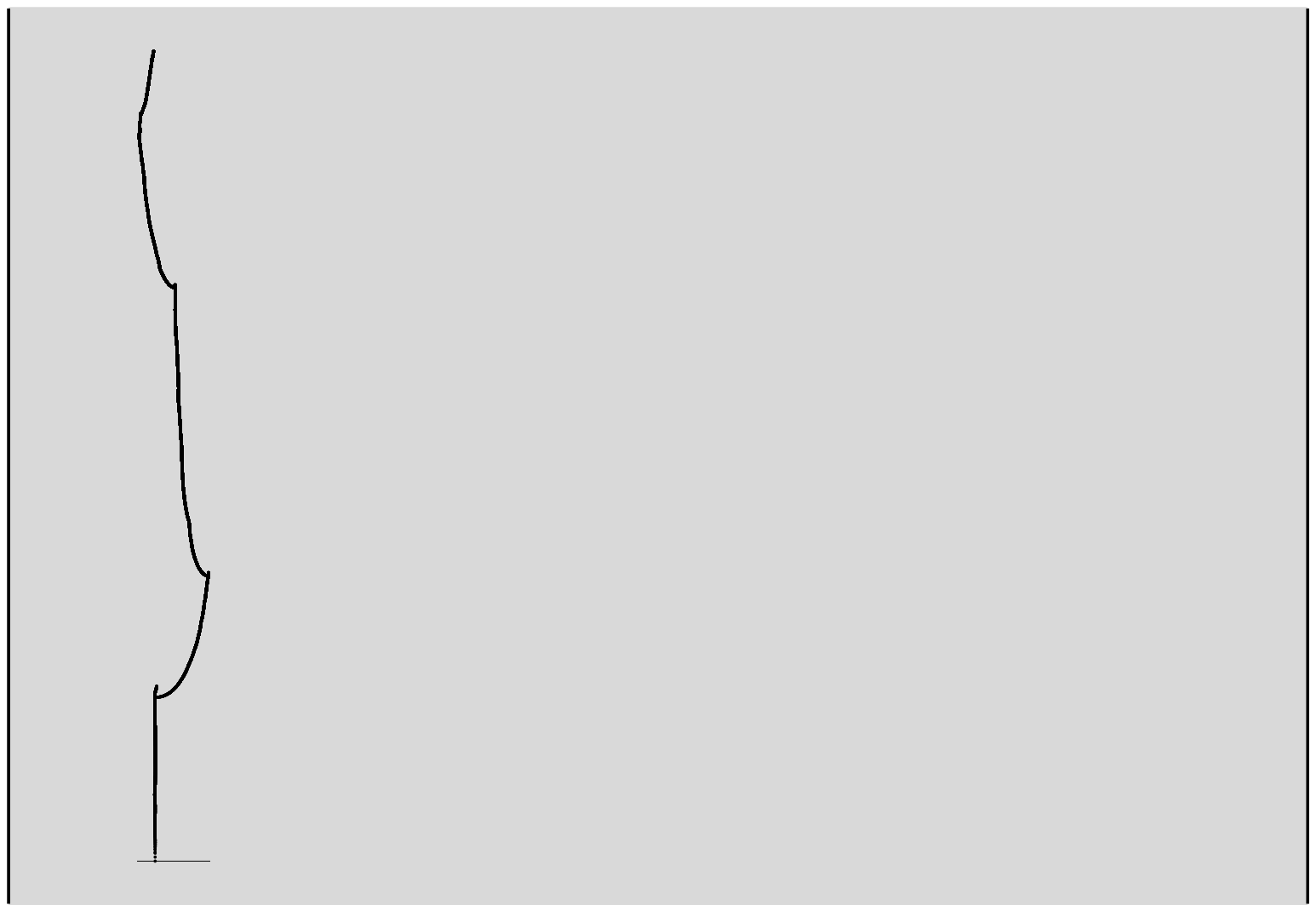}}

\centerline{Figure 1.1  ($\alpha=1$ and $t=0.1$)} }
\bigskip

\bigskip \vbox{ \epsfxsize=3.0in
  \centerline{\epsffile{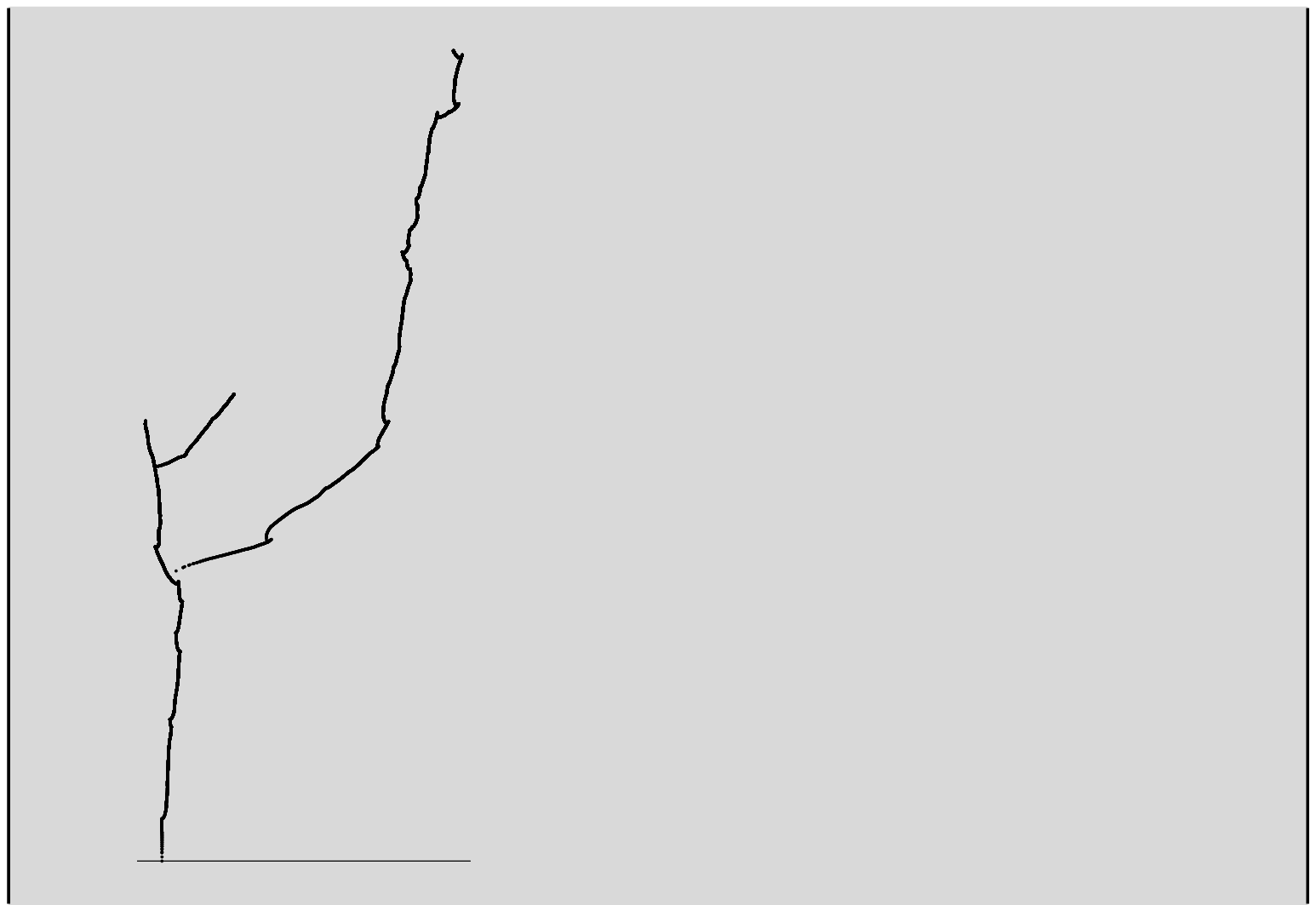}}

\centerline{Figure 1.2 ($\alpha=1$ and $t=1$) } }
\bigskip

\bigskip \vbox{ \epsfxsize=3.0in
  \centerline{\epsffile{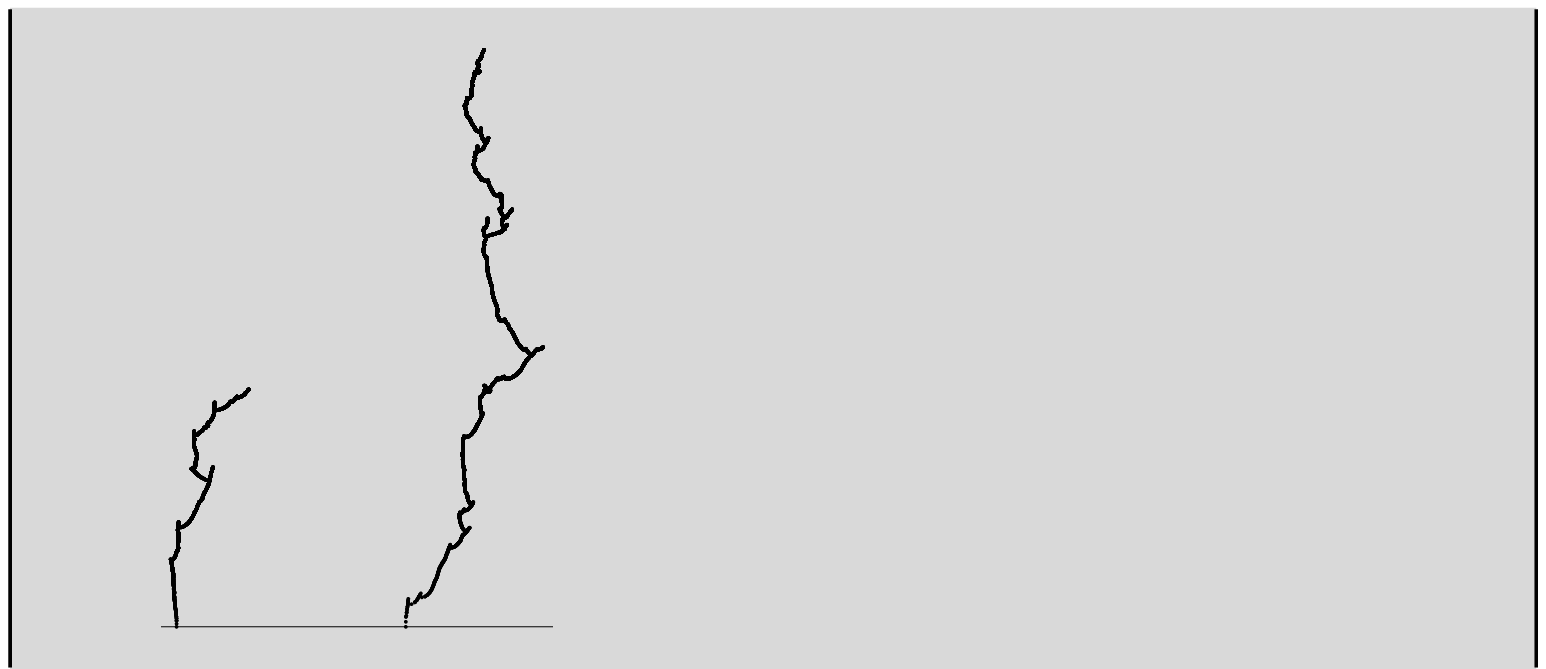}}

\centerline{Figure 1.3 ($\alpha=1$ and $t=10$) } }
\bigskip

\bigskip \vbox{ \epsfxsize=3.0in
  \centerline{\epsffile{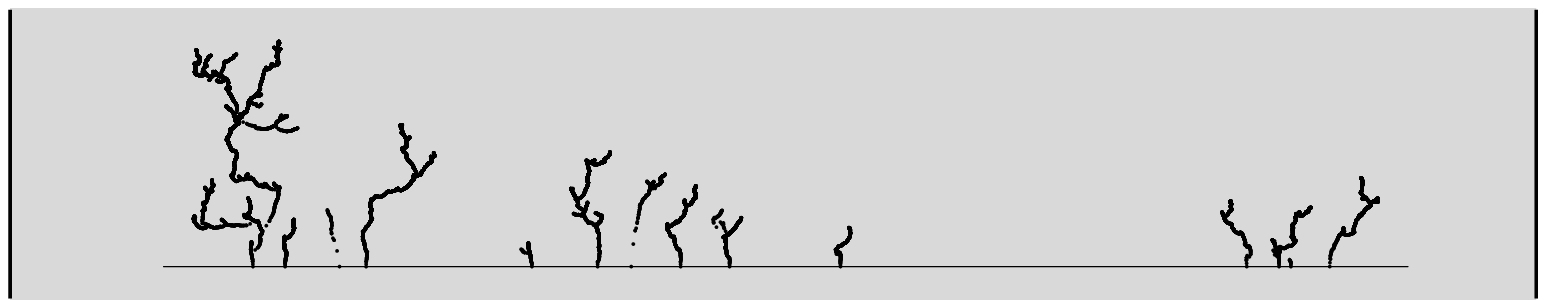}}

\centerline{Figure 1.4 ($\alpha=1$ and $t=100$) } }
\bigskip

As in the study of SLE in \cite{RS}, a key role in understanding
$K_t$ is therefore played by the derivative expectation
$\E[|f_t'(z)|^p]$. However in contrast with Brownian motion, the
infinitesimal generator of the symmetric $\alpha$-stable process $S$
on $\R$ is the fractional Laplacian
$\Delta^{\alpha/2}$, which is not very amenable
to calculations.
Many  nice smooth functions such as polynomials of
order 2 and beyond are not in its domain. For this technical reason,
we use the truncated symmetric standard $\alpha$-stable process $\wh
S$ instead, which is the symmetric $\alpha$-stable process $S$  with
jumps of size larger than $1$ removed. Any $C^2$-smooth function on
$\R$ is in the domain of the infinitesimal generator of $\wh S$.
Note that for the symmetric $\alpha$-stable process $S$, jumps of size
larger than 1 arrive according to a Poisson process. So there are
only a finite number of jumps of size larger than 1 in any given time
interval. For any $\kappa>0$,
 information on SLE driven by $S=\{S_t, t\geq 0\}$ can be easily
deduced from SLE driven by $\{S_{\kappa t}, t\geq 0\}$ (see Lemma
\ref{L:3.1} below), which in turn can be recovered from SLE driven
by  $\{\wh S_{\kappa t}, t\geq 0\}$ (see Lemma \ref{composition}
below).
Our main estimate here is
Theorem \ref{T:2.8}.
For $\kappa>0$, let  $W_t=\wh S_{\kappa t}$ and write
$$ f_t(z)-W_t=X_t+ i Y_t, \qquad t\geq 0.
$$
After a time change $\gamma_u:=\inf\{t\geq 0: Y_t\geq Y_0e^u\}$ and
$\wt f_u (z):=f_{\gamma_u}(z)$
we show in Section 4 that for every $0<p<2$ and $\delta>0$ there is
$\kappa>0$ such that for $W_t = \wh S_{\kappa t}$ and every $0<y<1$
$$
\E\left[|\wt f_{-\log y}'(z)|^{p}; \gamma_{-\log y} <\infty
\right]\leq C_{p,\delta} y^{-\delta}.
$$
In particular, this
  implies trivial integral means,
$$
    \beta(1)=0\quad {\text a.s.}
$$
for all $\kappa>0,$ and also for the (non-truncated) stable process.

\medskip

\noindent
 We apply the above derivative estimates to prove  in Section 5 that
 for every $T>0$, the maps of the
 backward flow  $f_t(z)$ of (\ref{eqn:1.2}) driven by $W_t=\wh S_{\kappa t}$ with
small $\kappa$ are
 uniformly $\gamma$-H\"older continuous on every bounded set $A\subset \H$
 for $t\in [0, T]$ with $\gamma$ close to
$1/6$.
The H\"older exponents are certainly not optimal
(we believe that the correct exponent is 1/2 for all $\alpha \in (0,
2)$). Nevertheless, this establishes enough regularity to prove that
the box counting (and hence the Hausdorff) dimension of the hull
$K_t$ (of SLE (\ref{eqn:1.1}) driven either by $W_t=S_{\kappa t}$ or
by $W_t=\wh S_{\kappa t}$ for every $\kappa >0$) is 1 a.s. It also
implies that the backward flow $f_t$ of (\ref{eqn:1.2}) driven
either by $W_t=S_{\kappa t}$ or by $W_t=\wh S_{\kappa t}$ for every
$\kappa >0$ is locally uniformly H\"older continuous in $\overline
\H$ a.s.
In particular, this implies that for each $t>0$,  the domain
$\H\setminus K_t$ is a H\"older domain almost surely.

\bigskip
\noindent
Finally, as another application of the H\"older
continuity of the maps of the backward flow $f_t(z)$ of
(\ref{eqn:1.2}),
 we prove that the trace is
right-continuous with left limits (RCLL in abbreviation):
 Let   $\{g_t, t\geq 0\}$ be SLE (\ref{eqn:1.1}) driven
either by $W_t=S_{\kappa t}$ or by $W_t= \wh S_{\kappa t}$.
We show in Theorem \ref{continuity} that for every $\alpha\in(0,2)$ and $\kappa
>0$, almost surely,
 for each $t>0$ the limit
$$\gamma(t) = \lim_{z\to W_{t}; z\in\H} g_t^{-1}(z)$$
exists, the function $t\mapsto \gamma(t)$ is RCLL, and $K_t =
\overline{ \gamma[0,t]}$.
This is achieved by first showing that with probability one, the
maps $\{g_t^{-1}, \, 0\leq t\leq T\}$ are equicontinuous on
$\overline \H$ for every $T>0$.

\bigskip
\noindent Independently from and parallel to this paper, Qing-Yang
Guan \cite{G} has recently investigated the continuity properties of
the trace of the Loewner equation driven by $W_t=B_{\kappa t} +
S_{\theta t}$ for $\kappa\geq 0,$ $\theta\geq 0,$ and $S$ the
symmetric $\alpha$-stable process with $0<\alpha<2$ (he informed us
that the assumption $\kappa>0$ in his manuscript is not needed).
Thus the main result of \cite{G} contains our Theorem
\ref{continuity} as the special case $\kappa=0.$  
 Whereas his proof
of the RCLL property is an adaptation of the continuity proof from
\cite{RS}, we employ a different simpler method that takes advantage
of the tree structure of the hulls (which works only for
$\kappa\leq4$), and is of independent interest.

\section{Definition and Basic Properties of symmetric $\a$-stable process }\label{s1}

A random variable $X$ is {\it symmetric $\a$-stable} if its
characteristic function

\begin{equation}\label{sas}
\E[ e^{i \theta X} ] = e^{- c|\theta|^{\alpha}}.
\end{equation}

\noindent For $\a=2,$ this is the normal distribution. It is not
hard to show (but nontrivial) that such $X$ exists if and only if
$0<\a\leq 2$ (see for instance \cite[Section 6.5]{Ch}).

\noindent
Write $X\sim S(\alpha,c)$ where $\alpha$ is called the {\it index}
and $c^{1/\alpha}$ is  the {\it scale}.
If $X_i\sim S(\alpha,c_i)$ are independent, then (\ref{sas}) immediately gives

$$X_1+X_2 \sim S(\alpha,c_1+c_2)$$
and
$$a X \sim S(\alpha, c a^{\a}).$$

\noindent The symmetric $\a$-stable process $S=\{S_t, t\geq 0\}$ (or
$\a$-stable L\'evy motion) is a L\'evy-process (meaning $S$ is right
continuous with left limits, and has stationary independent
increments) with $S_t-S_r$ are distributed according to the {\it
$\a$-stable law}:
$$S_t-S_r \sim S(\alpha,t-r)$$
for $0\leq r\leq t.$ Notice that stable process is self-similar: for
every $c>0$,

\begin{equation}\label{scaling}
\left\{ S_{c t} -S_0; \, t\geq 0 \right\} \doteq \left\{
c^{1/\alpha}\left( S_t-S_0\right); \, t\geq 0 \right\}
\end{equation}
where $\doteq$ denotes equality in distribution. This is the analog
of the classical Brownian scaling for Brownian motion. The {\it
transition density function} can be obtained from the characteristic
function by the inverse Fourier transform:
$$p(t,x, y)=p(t, x-y) = \P_x\left( S_t\in[y,y+dy] \right) /dy =
\frac1{2\pi}\int_{\R} e^{-i (x-y\theta} e^{- t|\theta|^{\alpha}}
d\theta.
$$
Explicit formulas for $p$ exist only in a few special cases (for
$\alpha=2$  we have the normal distribution
$p(t,x)=\exp(-x^2/2t)/\sqrt{2\pi t},$ and for $\alpha=1$ the Cauchy
distribution $p(t,x)=\frac{t}{\pi(t^2+x^2)}$). However we have the
following estimate (see, for example, \cite{CK})
\begin{equation}\label{pestimate}
p(t,x, y) \asymp t^{-\frac1{\a}} \wedge \frac{t}{|x-y|^{1+\a}}=
t^{-\frac1{\a}} \left(1\wedge\frac{ t^{\frac1{\a}} }{ |x-y|}
\right)^{1+\alpha}, \qquad t>0, \, x, y\in \R .
\end{equation}
Here for $a, b\in \R$, $a\wedge b:=\min \{a, b\}$  and $a \vee b =
\max \{a, b\}$. Therefore
\begin{equation}\label{pestimate2}
\P\left( |S_t-S_0|\geq x \right) \ \asymp \ 1\wedge
\frac{t}{|x|^{\a}}.
\end{equation}

\noindent We thus see that for $\a<2$, $S_t$ has infinity  variance,
and for $\a\leq 1$, $S_t$  is even does not integrable.

\begin{lemma}\label{max}
For $t>0$ and $x>0$,
$$\P \left( \max_{0\leq r\leq t} S_r > x \right) \leq 2 \P\left(S_t >x \right)
 \ \asymp\  1\wedge  \frac{t}{|x|^{\a}}.$$
\end{lemma}

\noindent {\bf Proof.} Let $T=\inf\{s: S_s>x\}$. Then $\P(T\leq t) =
\P \left( T\leq t, S_t\leq S_T \right) +\P \left( T\leq t, S_t>
S_T\right)\leq 2 \P \left( T\leq t, S_t\geq  S_T\right)\leq 2 \P
\left( S_t
>x \right).$ \karo

\bigskip

\noindent
For a Borel measurable function $f$ on $\R$, we define the
fractional Laplacian $\Delta^{\alpha/2} f=-(-\Delta)^{\alpha/2}f$ at
$x\in \R$ as follows:
$$ \Delta^{\alpha/2} f(x):= c_\alpha \, \lim_{\eps\downarrow 0}
\int_{\{|h|>\eps\}} \frac{f(x+h)-f(x)}{|h|^{1+\a}} dh,
$$
whenever the limit exists. It is easy to see that for  every $f\in
C^2_b (\R)$ and every $\delta >0$,
$$
\Delta^{\a/2} f (x)=   c_\alpha \, \int_{\R} \frac{f(x+h)-f(x)-f'(x)
\1_{\{|h|\leq \delta\}}}{|h|^{1+\a}} dh,
$$
which is well-defined and is in fact a bounded continuous function
in $x$. By Ito's formula
 \begin{equation}\label{eqn:2.6}
 t\mapsto
f(S_t)-\int_0^t \Delta^{\alpha/2}f(S_r) dr
\end{equation}
is a martingale for every $f\in C^2_b(\R)$ (see, e.g., the proof of
Proposition 4.1 in \cite{BC} for details). Let ${\cal A}$ be  the
Feller generator (that is, the infinitesimal generator in the space
$C_b(\R)$ of  bounded continuous functions  equipped with the supremum
norm $\| \cdot \|_\infty$) of the symmetric $\alpha$-stable process
$S$. Then the above implies that for $f\in C^2_b(\R)$,
\begin{equation}\label{generator}
{\cal A} f(x) = \lim_{t\to0}\frac{\E[f(x+S_t)]-f(x)}{t} =
\Delta^{\a/2}f .
\end{equation}
That is, $C^2_b(\R)\subset {\cal D}({\cal A})$ and for $f\in
C^2_b(\R)$, ${\cal A} f=\Delta^{\alpha/2} f$.

\medskip

 \noindent Let the domain $D\subset \R$ and let $f$ be
defined in all of $\R$ and continuous in $D$. Then $f$ is {\it
harmonic in $D$ with respect to $S$} if $f$ has the mean value
property
$$f(x) = \E_x[f(S_{\t_{B(x,r)}})]$$
for all balls $B(x, r)$ with closure in $D$, where $\tau_{B(x,
r)}:=\inf \left\{t\geq 0: S_t\notin B(x, r)\right\}$. Then the ball
can be replaced by any open $D_1$ with $\overline D_1\subset D$, see
\cite[Theorem 2.2]{CS}.

\vskip.5in
\noindent
The function
$$u(x) =
\begin{cases}\ |x|^{\a-1} \ \ \ \ \ &\text{if } \a\neq 1\\
\ \log |x| &\text{if } \a=1
\end{cases} $$
is harmonic in $\R\setminus\{0\}$ as is shown in \cite{ROKG} (for
$\a\neq1$ this follows from the harmonicity of the Kelvin transform
   $|x|^{\a-1} h(1/x)$ of the constant function $h\equiv1$, cf. \cite{L}).

\vskip.5in
\noindent
This can be used to obtain a quick proof of the recurrence resp. transience of
the stable process for $\a>1$ resp. $\a<1:$

From Ito's formula, $\{u(S_t), t\in [0, T_0)\}$ is a non-negative
local martingale, where $T_0=\inf\{t\geq 0: S_t=0\}$ and by Fatou's
Lemma it is a supermartingale. For $0\leq r< S_0=x <R$, we therefore
have
\begin{equation}\label{dynkin}
u(x)\geq \E_x \left[ u(S_{ {\t_r}\wedge{\t_R} }) \right] = \E_x
\left[u(S_{\t_R}) \1_{\{\t_R<\t_r\}} + u(S_{\t_r})
\1_{\{\t_r<\t_R\}} \right]
\end{equation}
where $\t_r=\inf\{t:|S_t|<r\}$ and $\t_R=\inf\{t:|S_t|>R\}$. For
$\a>1$ we get
$$\P_x (\t_R < \t_r ) \leq \frac{u(x)}{u(R)} \to 0 \quad{\text as}\quad R\to \infty$$
proving recurrence. Whereas for $\a<1$ we have
$\P_x (\t_r < \t_R ) \leq \frac{u(x)}{u(r)}$ and so after letting
$R\to \infty$, we get
$$ \P_x (\t_r < \infty) \leq \frac{u(x)}{u(r)}<1,
$$
giving the transience of $S$.
Moreover, if we let $r\downarrow 0$
in the last formula,
we have  for $\a <1$ and every $x\not=0$
$$ \P_x (\sigma_{\{0\}} <\infty) =0.
$$
Here
$\sigma_{\{0\}}=\inf\{t>0: S_t=0 \hbox{ or } S_{t-} =0\}$.
In other words, almost surely neither $S_t$ nor $S_{t-}$ will
visit 0.

\section{The Loewner equation}\label{s2}

\subsection{Deterministic equation}\label{S:3.1}

Let $W_t$ be a real-valued function that is right continuous with left limits, RCLL for short.
For each initial point $z\in\C\setminus\{0\}$,
the {\it Loewner differential equation}

\begin{equation}\label{LE}
\partial_t g_t(z) = \frac2{g_t(z)-W_t} \ ,\quad g_0(z)=z
\end{equation}
has a unique solution up to a time $0<T_z\leq\infty$ where
$g_t(z)=W_{t-}$ or $g_t(z)=W_t$.
  More precisely, let
  $$ T_z=\sup\left\{t: \inf_{s\in[0,t]}|g_s(z)-W_s|>0 \right\},
  $$
   then the initial value problem
(\ref{LE}) has a unique solution on $[0,T_z)$ and if $T_z<\infty$
then either
 $\displaystyle \liminf_{t\to T_z -} |g_t(z)-W_t| = 0,$ or
$\displaystyle \liminf_{t\to T_z-}
  |g_t(z)-W_t| >0$ and $g_{T_z}(z)=W_{T_z}$ (in this case, $W$ jumps
at time $T_z$). The subset
$$K_t = \{z\in\overline\H: T_z\leq t\}$$
is a compact subset of the closed upper half plane $\overline\H$ and
is called  the {\it hull} of LE (\ref{LE}).  It is well-known that
the map $z\mapsto g_t(z)$ is a conformal map (i.e. analytic and
one-to-one) from $\H\setminus K_t$ onto $\H,$ with Laurent series
$g_t(z) = z + \frac{2t}{z} + O(\frac1{z^2})$ near $\infty.$ From the
uniqueness of normalized conformal maps it follows that
$K_t\cap\H\neq\emptyset$ is strictly increasing in $t.$ Writing
$g_t(z) = x_t+i y_t$ and taking real- and imaginary parts in
(\ref{LE}), the Loewner equation reads
\begin{equation}\label{realLE}
\begin{aligned} \partial_t x_t &=  &2\frac{x_t-W_t}{(x_t-W_t)^2 + y_t^2} \\
\partial_t y_t &= &-2\frac{y_t}{(x_t-W_t)^2 + y_t^2}
\end{aligned}
\end{equation}

\noindent
It is easy to see that when $W\equiv 0$, $g_t(z)=\sqrt{z^2+4t}$ and
$K_t=\gamma [0, t]$, where $\gamma (t)= i\, 2\sqrt{t} $.

\subsection{LE driven by stable processes}

\noindent
As $t\to\infty,$ the diameter of the hulls $K_t$ tends to infinity.
In fact,
\begin{equation}\label{diam}
\frac1C \sqrt t \leq \diam K_t \leq C (\sqrt t + \sup_{0\leq r\leq s\leq t}|W_r-W_s|)
\end{equation}
for some universal $C$ and all $t.$ What do the
 hulls of LE driven by stable process (that is, $W$ in (\ref{LE})
 is a symmetric $\alpha$-stable process)
look like if we scale them back down as $t\to\infty$, or scale them
up when $t\to0$? We will see that both the ``conformally natural"
and the "metrically natural" way of rescaling the hulls does not
lead to any interesting sets: If we scale them so as to have
(halfplane-) capacity one or so that the diameter is one, then the
hulls converge to a vertical line segment as $t\to0$ and to the
empty set as $t\to\infty.$ To make this precise, let  $c>0.$ The
solution to (\ref{LE}) with
\begin{equation}\label{bs1}
\wt W_t = \frac1c W_{c^2 t}
\end{equation}
is given by the function $\wt g_t(z) = \frac1c g_{c^2 t}(c z).$ It
follows that the hulls are related by

\begin{equation}\label{bs2}
\wt K_t = \frac1c K_{c^2t}.
\end{equation}
If $W_t$ is a Brownian motion with variance $\kappa$, then $\frac1c
W(c^2 t)$ has the same distribution which translates to the
important and useful scaling invariance of the $SLE$ hulls. If $W_t$
is $\a$-stable then $\frac1c W(c^2 t)$ is $\a$-stable too but the
scale is different: From (\ref{scaling}) it follows that

\begin{equation}\label{bs3}
\frac1c W(c^2 t)\doteq c^{\frac2{\a} -1}W_t.
\end{equation}

\noindent
Let $\{g_t(z), t\geq 0\}$ be the SLE driven by  $W=S$, the symmetric
standard $\alpha$-stable process on $\R$, with hulls $\{K_t, t\geq 0\}$.
For $c>0$, define $\wh g_t(z):=c^{-1} g_{c^2t} (cz)$. Then
$$ \partial_t \wh g_t(z)=\frac{2}{\wh g_t(z)-c^{-1}S_{c^2t}} \qquad
\hbox{with } \wh g_t(z)=z.
$$
So  $\{\wh g_t(z), t\geq 0\}$ is SLE with hulls $\{\wh K_t=c^{-1}K_{c^2 t}, t\geq 0\}$,
 driven by symmetric $\alpha$-stable process $\{c^{-1}S_{c^2t}, t\geq 0\}
\doteq \{ S_{c^{2-\alpha}t}, t\geq 0\}$ running at a different speed.
We record this as a lemma for future reference.

\begin{lemma}\label{L:3.1} Let $\{K_t, t\geq 0\}$ be the hulls of SLE driven by
$W=S$, the symmetric standard $\alpha$-stable process. Then for every $c>0$,
$\{c^{-1}K_{c^2 t}, t\geq 0\}$ has the same distribution as the hulls of SLE
driven by $\{W_t=S_{c^{2-\alpha}t}, t\geq 0\}$. Hence the geometric information on
hulls of SLE driven by
$W_t=S_t$ and by $W_t=S_{\kappa t}$ can be deduced one from the other.
\end{lemma}

\noindent
From (\ref{bs3}), it is not difficult to prove:

\begin{prop}\label{p1}
Let $0<\a <2$ and $\{K_t, t\geq 0\}$ be the hulls of SLE driven by $W=S$.
As $s \to0,$ the rescaled hulls $\frac1{s} K_{s^2}$ converge to the
vertical line segment $[0,2i]$ (in the Hausdorff metric) in
probability. On the other hand, for all $\eps>0$,
$$
\lim_{s\to \infty} \P\left( \frac1s K_{s^2} \cap \{y>\eps\} \not=
\emptyset \right) =0.
$$
\end{prop}

\noindent
The proof uses the following simple result for deterministic hulls:

\begin{lemma}\label{l1}

\noindent
(a) If $W_t\in[a,b]$ for all $t\in[0,T]$, then $K_T\subset [a,b]\times \R.$

\noindent (b) Let $0<\eps<1$ and $r>1.$ If $I\subset\R$ is an
interval of length $\sqrt {T}$ and $10 I$ the concentric interval of
size $10\sqrt{T}$, and if $$\int_0^T \1_{\{W_t\in 10I\} }dt \leq
\eps T,$$ then
$$K_T\cap \ I\times[4\sqrt{\eps T},\infty) = \emptyset.$$
\end{lemma}

\noindent {\bf Proof.} (a) If $z=x+iy\in\H$ with $x<a$ (resp. $>b$),
then $\partial_t x_t(z) < 0$ (resp. $>0$) and hence $|g_t(z)-W_t|$
is bounded from below by $|x-a|$.

\noindent
(b) By means of Brownian scaling (\ref{bs1}) and (\ref{bs2}), we may
assume $T=1$ and $I=[-1/2,1/2]$. Fix $z_0\in
I\times[4\sqrt{\eps},\infty)$ and write $g_t(z_0)=x_t+i y_t.$ We may
assume $y_0<2$, else trivially $z_0\notin K_1$ (only the hull of the
constant function $W_t\equiv 0$ reaches height 2). Let $T_1\leq1$ be
maximal time such that $x_t+iy_t \in [-2,2] \times
[2\sqrt{\eps},\infty)$ for all $t\in[0,T_1]$. We will show $T_1=1$
and hence $z_0\notin K_1$, proving the lemma. Up to $T_1$, from
(\ref{realLE}) we have
\begin{equation}
\begin{aligned}
|x_t-x_0| &\leq \int_0^{T_1} 2\frac{|x_t-W_t|}{(x_t-W_t)^2 + y_t^2} dt \\
& \leq \int_{\{W_t\in 10 I\}}2\frac1{2y_t}dt + \int_{\{W_t\notin 10
I\}}2\frac1{|x_t-W_t|}dt \\
& \leq \frac{\eps}{4\sqrt{\eps}} + \frac23<\frac32
\end{aligned}
\end{equation}
Thus $x_t$ does not reach the boundary of $[-2,2]$. Similarly,

\begin{equation}
\begin{aligned}
|y_{T_1} - y_0| &= \int_0^{T_1} 2\frac{y_t}{(x_t-W_t)^2 + y_t^2} dt \\
&\leq \int_{\{W_t\in 10 I\}}\frac2{y_t}dt + \int_{\{W_t\notin 10 I\}}2\frac{y_t}{(x_t-W_t)^2}dt\\
&\leq \eps\frac2{y_0} + 2\frac{y_0}{9}\leq \sqrt{\eps} +
\frac{y_0}{4}<\frac{y_0}{2}
\end{aligned}
\end{equation}
and therefore $y_t$ does not reach $2\sqrt{\eps}$. Hence $T_1=1$ and
the lemma is proved. \karo

\vskip.5in \noindent {\bf Proof of Proposition {\ref{p1}}}. From
(\ref{bs3}) we have that the rescaled hull has the same distribution
as the time 1 hull of the map $t\mapsto s^{\frac2{\a} -1}W_t.$ By
Lemma {\ref{max}}, the support of this function tends to zero in
probability, and from Lemma \ref{l1} (a) it follows that the width
of the hull tends to zero in probability. Since the
halfplane-capacity it 1, the height has to converge to 2 and the
hull converges to the segment $[0,2i]$ as $s\to0.$ By Lemma \ref{l1}
(b), the second claim is equivalent to saying that the maximal
amount of time that $W_t$ spends in an interval $[x,x+\delta]$ tends
to zero in probability as $\delta\to0.$ \karo

\subsection{Stable LE on $\R$}
\noindent
When $z$ is a non-zero real number, $Z_t:=g_t(z)$ of (\ref{LE}) is
real-valued.
We will call the real-valued equation
\begin{equation}\label{LE+}
\partial_t Z_t = \frac2{Z_t-W_t}
\end{equation}
the {\it forward Loewner equation on $\R$ driven by } $W_t,$
 and
\begin{equation}\label{LE-}
\partial_t Z_t = -\frac2{Z_t-W_t}
\end{equation}
 the {\it backward Loewner equation on $\R$}.
 The latter corresponds to the backward flow $f_t(z)$ of
 (\ref{eqn:1.2}) with $z\in \R\setminus \{0\}$.
 If $W$ is the symmetric $\a$-stable process on $\R$,
 then
the generator of $$X_t=Z_t-W_t$$ in the forward resp. backward
equation is
$$A_{\pm} = \pm \frac2{x} \frac{d}{dx} - (-\Delta)^{\a/2}.$$
For the $(-\Delta)^{\a/2}$-harmonic function $u(x) =|x|^{\a-1}$ ($\a\neq1$)
we have
$$A_{\pm} u = \pm 2(\a-1)|x|^{\a-3}.
$$
Thus $u$ is superharmonic for $A_+$ if $\a<1$ and for $A_-$ if $\a>1.$
With the above reasoning (\ref{dynkin}) we obtain
\begin{prop}\label{p2}
For $\a>1$, $X_t$ is recurrent in the
backward LE on $\R$, whereas for $\a<1,$ $X_t$ is transient in the
forward LE on $\R$
and almost surely, neither $X_t$ nor $X_{t-}$ visits $0$.
\end{prop}

\noindent
Notice that in
 SLE driven by Brownian motion,
which corresponds to the case $\a=2$,
 the question of recurrence versus transience of $X_t$ in the forward LE is rather subtle:
 If $B_t$ is Brownian motion and
$W_t=\sqrt{\kappa} B_t$ with $\kappa\leq4,$ we have transience whereas for $\kappa>4$
 we have recurrence.

\bigskip
\noindent
We will now prove a partial converse to Lemma \ref{l1}(a)
about the deterministic forward LE (\ref{LE}) in $\H$.
If a point $x_0$ on the real line stays away by $\eps$ from the
singularity, then the disc of radius $\eps$ at this point does not
meet the singularity and therefore is disjoint from the hull. More
generally, if the {\it real part} of some point stays away by
$\eps$, then the $\eps$- disc around this point is disjoint from the
hull.

\bigskip

\noindent Let $g_t(z)$ be the solution to the deterministic LE
(\ref{LE})  with $z\in \overline \H\setminus \{ 0\}$ and define
$X^z_t=g_t(z)-W_t$. When $z\in \R\setminus \{0\}$,
$Z_t = g_t(z)$ solves the forward LE (\ref{LE+}) on $\R$
discussed at the beginning of the section and
$X^z_t=Z_t-W_t$ is real-valued. When $z\in \H$, $X^z_t$ is complex
valued.
 We will use
$B(z, r)$ to denote the ball in $\R^2=\C$) centered at $z$ with
radius $r$.

\begin{lemma}\label{l2}
  If $|\Re X^{z_0}_t|\geq \eps$ for some $z_0\in\overline\H$
and all $0\leq t \leq T,$ then $|\Re X^{z}_t| > 0 $ for all $z\in
B(z_0, \eps) \cap \overline \H$ and all $0\leq t \leq T.$ In
particular, $B(z_0,\eps)\cap K_T = \emptyset.$
\end{lemma}

\noindent {\bf Proof.} From (\ref{LE}) we have
$$\partial_t (X^{z}_t - X^{z_0}_t ) =
\frac{2}{X^z_t} -\frac{2}{X^{z_0}}=2\frac {X^{z_0}_t -
X^{z}_t}{X^{z}_t X^{z_0}_t}
$$
 and so
\begin{equation}
\begin{aligned}
\partial_t |X^{z}_t - X^{z_0}_t|^2  &=
2 \Re \left[ \partial_t (X^{z}_t - X^{z_0}_t ) \overline{(X^{z}_t - X^{z_0}_t )} \right] \\
&= -4\frac{|X^{z}_t - X^{z_0}_t|^2}{|X^{z}_t X^{z_0}_t|^2} \Re
(X^{z}_t X^{z_0}_t).
\end{aligned}
\end{equation}
It follows that $|X^{z}_t - X^{z_0}_t|$ is decreasing  because $\Re
(X^{z}_t X^{z_0}_t)>0$ as long as $|\Re X^{z_0}_t|\geq \eps$ and
$|\Re (X^{z_0}_t-X^{z}_t)| < \eps$. Since $|\Re X^{z_0}_t|\geq \eps$
for every $0\leq t \leq T,$ we have for every $z\in B(z_0, \eps
)\cap \overline \H$,
$$|X^{z}_t - X^{z_0}_t|\leq |z-z_0|<\eps \qquad \hbox{for every }
0\leq t\leq T
$$
and so $|\Re X^{z}_t| > 0 $ for every  $0\leq t \leq T.$  \karo

\bigskip

\noindent
Now suppose $x\in\R\setminus\{0\}$ and $g_t(x)$ is the solution to the LE (\ref{LE})
driven by a symmetric $\a$-stable process $W$ on $\R$ with $\a<1$.
As mentioned previously, $g_t(x)$ is the solution to the forward LE on $\R$.
Proposition \ref{p2} tells us that for $X^x_t=g_t(x)-W_t$,
$$
  r:=\inf_{t\geq 0} |\Re X_t^x|=\inf_{t\geq 0} |X_t^x|>0 \qquad
  \hbox{a.s.}
$$
 We then have by Lemma \ref{l2}
$$B(x,r)\cap \bigcup_{t>0} K_t = \emptyset \qquad \hbox{a.s.} $$

 \bigskip

\section{Derivative Estimates}\label{Derivative Estimates}

\noindent
We would like to estimate the derivative of $h_t=g_t^{-1}.$ Because
$h_t$ satisfies the PDE
$$\partial_t h_t(z) = -2\partial_z h_t(z)/(h(t(z)-W_t)
$$
 rather than an ODE, it is usually easier to work with the time $t$
map $f_t$ of the backward Loewner equation
 (\ref{eqn:1.2}).
The connection is as follows: If $g_t$ is the solution to (\ref{LE})
driven by a function $W_t$ ($0\leq t\leq T$), and if $f_s$ is the
solution to (\ref{LE-}) driven by $\wt W_s= W_{T-s},$ then $f_T =
g_T^{-1}.$ But generally $f_t\neq g_t^{-1}$ for $t<T.$ Because for
the symmetric stable process, $s\mapsto W_{T-s}-W_T$ has the same
distribution as $W_s,$ it follows that for each fixed $T>0$, the
random conformal map $f_T(z)$ of $\H$ has the same distribution as
$g_T^{-1}(z-W_T)+W_T$ (but the family of maps, $\{f_t(\cdot ), t\geq
0\}$ does not have the same distribution as $\{g^{-1}_t(\cdot
-W_t)+W_t, t\geq 0\}$). For the remainder of this section, we
consider the time $t$ map $f_t$ of the backward Loewner equation
 (\ref{eqn:1.2}).

 \medskip

Let $(X_t, Y_t):=Z_t-W_t$. Then by (\ref{LE-}),
$$ d (X_t+i Y_t)= \frac{-2}{X_t+iY_t} dt -dW_t=\frac{-2X_t+2iY_t}{X_t^2+Y_t^2} dt -dW_t.
$$
Hence
\begin{equation}\label{eqn:2.10}
dX_t=-\frac{2X_t}{X_t^2+Y_t^2}dt -dW_t \quad \hbox{ and } \quad
  dY_t=\frac{2Y_t}{X_t^2+Y_t^2} dt.
  \end{equation}
  In particular, we have  $d\ln Y_t= \frac2{X_t^2+Y_t^2} dt$ and
  so
  $$ Y_t=Y_0 e^{\int_0^t \frac2{X_t^2+Y_t^2} dt}.
  $$
We record a simple lemma for later use. Let
$\phi_t(z)=\sqrt{z^2-4t}$ be the solution to the backward LE
(\ref{eqn:1.2}) driven by the constant function $W\equiv0.$

\begin{lemma}\label{height}
For every $Z_0=X_0+iY_0$ with $Y_0\in (0, 1]$,
$$ Y_t \leq \Im \phi_t(i Y_0) \leq
 \sqrt{1+4t}
\qquad \hbox{for every } t>0.
$$
\end{lemma}

\noindent
\pf From (\ref{eqn:2.10}) we have $d Y_t \leq 2dt/Y_t$ with equality if and only if
$X_t\equiv 0$ and therefore $W_t\equiv 0.$ Thus $d(Y_t^2)\leq 4dt$ and integration gives
$Y_t^2\leq Y_0^2 + 4t$ with equality if and only if $W_t\equiv 0.$

\qed

\medskip

\noindent
For $u>0$, define
 \begin{equation}\label{gamma}
\gamma_u =\inf \left\{t>0:
Y_t\geq Y_0 e^u\right\} = \inf \left\{t>0: \int_0^t \frac2{X^2_s+Y^2_s} ds
\geq u\right\}.
\end{equation}

\begin{thm}\label{L:2.5} Let $W_t=S_t$  be a standard symmetric
 $\alpha$-stable process on $\R$
{\rm (}that is, $W\sim S(\alpha, 1))$.
 Then for  every  $z=x+iy\in \H$ and $u>0$, $\P_z(\gamma_u <\infty) =1$
when $\alpha \in [1, 2)$ and $\P_z(\gamma_u=\infty)>0$ when $\alpha
\in (0, 1)$.
\end{thm}

\pf Define $u_0:=\inf\{u: \, \gamma_u=\infty\}$, and $(\wt X_u, \wt
Y_u):=(X_{\gamma_u}, Y_{\gamma_u})$ for $u<u_0$. Clearly for
$u<u_0$, $\wt Y_u=Y_0e^u$. Note that under $\P_z$, $(X_0, Y_0)=(x,
y)$, so for $u<u_0$, $\wt Y_u=ye^u$ and
 \begin{equation}\label{Xu}
\wt X_u=X_{\gamma_u}=x-\int_0^{\gamma_u} \frac{2X_s}{X^2_s+Y^2_s} ds
-W_{\gamma_u} = x-\int_0^u \wt X_s \, ds -W_{\gamma_u}.
 \end{equation}
 By  \cite[Theorem 3.1]{RW}, there is a symmetric $\alpha$-stable
process $Z$ on $\R$ such that
$$   W_{\gamma_u} = \int_0^u \left( \frac{\wt X_r^2+
y^2 e^{2r}}{2}\right)^{1/\alpha}  dZ_r \qquad \hbox{on } [0, u_0).
$$
 Thus $\wt X$
satisfies the following SDE
 \begin{equation}\label{eqn:2.11}
  d\wt X_u=  - \wt X_u du -\left( \frac{\wt X_u^2+ y^2
e^{2u}}{2}\right)^{1/\alpha}  dZ_u \qquad \hbox{on } [0, u_0) \qquad
\hbox{with } \wt X_0=x,
\end{equation}
 where $Z$ is a symmetric $\alpha$-stable
process on $\R$. We can rewrite (\ref{eqn:2.11}) as
$$ d(e^u \wt X_u)=- e^{(1-2/\alpha)u}\left( \frac{(e^u  \wt X_u)^2+ y^2
e^{4u}}{2}\right)^{1/\alpha} dZ_u .
$$
By \cite[Lemma 4.5 and Theorem 4.6]{EK}, the above SDE for
$U_t:=e^t\wt X_t$ has a unique weak solution. Moreover
\cite[Theorems 4.7 and 4.9]{EK} tell us that the solution has
non-explosion if and only if $\alpha \in [1, 2)$ (see also \cite{PZ}
for the case of $\alpha \in (1, 2)$). It follows that SDE
(\ref{eqn:2.11}) has a unique weak solution $\overline X$ that has
infinite lifetime if and only if $\alpha \in [1, 2)$. Note that the
process $(\overline X, ye^u)$ extends $(\wt X_u, \wt Y_u)$ in law.
So we have for $\alpha \in [1, 2)$, $u_0=\infty$ a.s., in other
words, for any $t>0$, the original height process $Y$ can reach
level $ye^{t}$ with probability 1. When $\alpha \in (0, 1)$, the
proof of \cite[Theorem 4.9]{EK} illustrates that
$\P_z(\gamma_t=\infty)>0$ for every $t>0$. This proves the lemma.
\qed

\bigskip
\noindent
As we mentioned in the Introduction, many
smooth functions such
as polynomials of order 2 and higher are not
 $\Delta^{\alpha/2}$-differentiable.
For this reason, we need to look at truncated symmetric stable processes.
Let
$$\wh S_t:=S_t-\sum_{0< r \leq t} (S_r-S_{r-}) \1_{\{|S_r-S_{r-}|>1\}}, \qquad t\geq
0.
$$
The process $\wh S$ is a L\'evy process with L\'evy characteristic
measure $c_\alpha |h|^{-1-\alpha} \1_{\{|h|\leq 1\}}$ (see
\cite[Theorem I.1]{Be}). We call $\wh S$ a {\it truncated symmetric
standard $\alpha$-stable process} with jumps of size larger than $1$
removed. Define for $f\in C^2(\R)$,
$$ \wh \Delta^{\alpha/2} f(x):=
\int_{-1}^1 (f(x+h)-f(x)-f'(x)h ) c_\alpha |h|^{-1-\alpha} dh.
$$
Note that by Taylor expansion, we have for $f\in C^2(\R)$,
\begin{equation}\label{eqn:estimate}
| \wh \Delta^{\alpha/2} f(x) | \leq \sup_{w\in [x-1, x+ 1]} |
f''(w)| \frac12 \int_{-1}^1 c_\alpha |h|^{1-\alpha} dh= C_0
\sup_{w\in [x-1, x+1]} | f''(w)|.
\end{equation}
Thus there is a richer family of test functions at our disposal for
$\wh \Delta^{\alpha/2}$  than for $\Delta^{\alpha/2}$. Just as in
(\ref{eqn:2.6}), it follows from Ito's formula that for every $f\in
C^2(\R)$,
\begin{equation}\label{eqn:4.6}
t\mapsto f(\wh S_t) -\int_0^t \wh \Delta^{\alpha/2}f (\wh S_r) dr
\end{equation}
is a local martingale.

\bigskip

\begin{lemma}\label{L:2.7} Let $f(x):=(x^2+a^2)^{p/2}$ where $a>0$ and $p<2$.
Then there are
constants $C_1, C_2>0$ depending on $p$ and $\alpha$ only such that
\begin{equation}\label{eqn:2.13}
 | \wh \Delta^{\alpha/2} f(x) | \leq C_1
(x^2+a^2)^{(p-\alpha)/2} +C_2.
\end{equation}
When $p=\alpha,$ the right-hand side is to be interpreted as $\log(1/(x^2+a^2).$
\end{lemma}

\pf  The proof is similar to that of Lemma 2.9 in \cite{GW}. Assume
first that $|x| \leq a$. Then $w:=x/a \in [-1, 1]$. When $0<a<1/2$,
we have
\begin{eqnarray*}
| \wh \Delta^{\alpha/2} f(x)| &=& \left| \lim_{\eps \to 0} c_\alpha
\int_{\{\eps<|h|\leq 1 \}}
 \frac{ \left( (x+h)^2+a^2 \right)^{p/2} -(x^2+a^2)^{p/2}}{|h|^{1+\alpha}}
 dh \right|\\
 &\leq & \lim_{\eps \to 0} c_\alpha \left| \int_{\{\eps<|h|\leq 1\}}   a^p
 \frac{ \left( (\frac{x}{a}+\frac{h}{a})^2+1 \right)^{p/2} -\left(\frac{x^2}{a^2}+1 \right)^{p/2}}
  {|h|^{1+\alpha}} dh \right| \\
&=& \lim_{\eps \to 0} c_\alpha \left| \int_{\{\eps/a <|t|\leq 1/a\}}  a^{p}
 \frac{ \left( (w+t)^2+1 \right)^{p/2} -\left( w^2+1 \right)^{p/2}}
  {a^{1+\alpha} |t|^{1+\alpha}}\,  a dt \right| \\
  &=& a^{p-\alpha}  \lim_{\eps \to 0} c_\alpha  \left| \left( \int_{\{\eps/a<|t|\leq 2\}}
  +\int_{\{2<|t|\leq 1/a\}} \right)
 \frac{ \left( (w+t)^2+1 \right)^{p/2} -\left( w^2+1 \right)^{p/2}}
  {|t|^{1+\alpha}}\,   dt \right| \\
  &\leq & a^{p-\alpha}  c_\alpha   \left| \int_{\{ |t|\leq 2\}}
 \frac{ \left( (w+t)^2+1 \right)^{p/2} -\left( w^2+1 \right)^{p/2}
 -pw(w^2+1)^{p/2-1}   t} {|t|^{1+\alpha}}\,   dt \right| \\
 && + a^{p-\alpha}  c_\alpha \left| \int_{\{2<|t|\leq 1/a\}}
 \frac{ \left( (w+t)^2+1 \right)^{p/2} -\left( w^2+1 \right)^{p/2}}
  {|t|^{1+\alpha}}\,   dt \right| \\
  &\leq & c_1 a^{p-\alpha} + a^{p-\alpha} c_2 \int_2^{1/a}
  t^{p-1-\alpha} dt \\
  &\leq & c_1 a^{p-\alpha} +c_3 \\
  &\leq &c_4 (x^2+a^2)^{(p-\alpha)/2} +c_3.
\end{eqnarray*}
When $a\geq 1/2$, by the same calculation as above we have
\begin{eqnarray*}
| \wh \Delta^{\alpha/2} f(x)|  &=& \lim_{\eps \to 0} c_\alpha \left|
\int_{\{\eps/a <|t|\leq 1/a\}}
 \frac{ a^{p} \left( (w+t)^2+1 \right)^{p/2} -\left( w^2+1 \right)^{p/2}}
  {a^{1+\alpha} |t|^{1+\alpha}}\,  a dt \right|  \\
  &\leq & a^{p-\alpha}  c_\alpha   \left| \int_{\{ |t|\leq 1/a\}}
 \frac{ \left( (w+t)^2+1 \right)^{p/2} -\left( w^2+1 \right)^{p/2}
 -pw(w^2+1)^{p-1}   t} {|t|^{1+\alpha}}\,   dt \right| \\
  &\leq & c_1 a^{p-\alpha} \\
  &\leq & c_4 (x^2+a^2)^{(p-\alpha)/2}.
\end{eqnarray*}
This proves (\ref{eqn:2.13}) for the case of $|x|\leq a$.

Now assume $|x|>a$. Then $u:=a/x\in (-1, 1)$. If $a<|x|\leq 1/2$, we
have
\begin{eqnarray*}
| \wh \Delta^{\alpha/2} f(x)|
 &\leq & \lim_{\eps \to 0} c_\alpha \left| \int_{\{\eps<|h|\leq 1\}} x^2
 \frac{ \left( (1+\frac{h}{x})^2+\frac{a^2}{x^2} \right)^{p/2} -
 \left(1+\frac{a^2}{x^2} \right)^{p/2}}
  {|h|^{1+\alpha}} dh \right| \\
&=& \lim_{\eps \to 0} c_\alpha \left| \int_{\{\eps/|x| <|t|\leq
1/|x|\}}
 \frac{ |x|^{p} \left( (1+t)^2+u^2 \right)^{p/2} -\left( 1+u^2 \right)^{p/2}}
  {|x|^{1+\alpha} |t|^{1+\alpha}}\,  |x| dt \right| \\
  &=& |x|^{p-\alpha}  \lim_{\eps \to 0} c_\alpha  \left|
 \left( \int_{\{\eps/|x|<|t|\leq 1/2\}}
  +\int_{\{ 1/2<|t|\leq 1/|x|\}} \right)
 \frac{ \left( (1+t)^2+u^2 \right)^{p/2} -\left( 1+u^2 \right)^{p/2}}
  {|t|^{1+\alpha}}\,   dt \right| \\
  &\leq & |x|^{p-\alpha}  c_\alpha   \left| \int_{\{ |t|\leq 1/2\}}
 \frac{ \left( (1+t)^2+u^2 \right)^{p/2} -\left( 1+u^2 \right)^{p/2}
 -p(1+u^2)^{p-1}   t} {|t|^{1+\alpha}}\,   dt \right| \\
 && + |x|^{p-\alpha}  c_\alpha \left| \int_{\{1/2<|t|\leq 1/|x|\}}
 \frac{ \left( (1+t)^2+u^2 \right)^{p/2} -\left( 1+u^2 \right)^{p/2}}
  {|t|^{1+\alpha}}\,   dt \right| \\
  &\leq & c_5 |x|^{p-\alpha} +c_7\\
 &\leq & c_8 (x^2+a^2)^{(p-\alpha)/2} +c_7.
\end{eqnarray*}
When $|x|>\max \{1/2, a\}$, we have from above
\begin{eqnarray*}
| \wh \Delta^{\alpha/2} f(x)|
 &\leq &  \lim_{\eps \to 0} c_\alpha \left| \int_{\{\eps/|x| <|t|\leq
1/|x|\}}
 \frac{ |x|^{p} \left( (1+t)^2+u^2 \right)^{p/2} -\left( 1+u^2 \right)^{p/2}}
  {|x|^{1+\alpha} |t|^{1+\alpha}}\,  a dt \right| \\
  &=& |x|^{p-\alpha}   c_\alpha   \left| \int_{\{ |t|\leq 1/|x|\}}
 \frac{ \left( (1+t)^2+u^2 \right)^{p/2} -\left( 1+u^2 \right)^{p/2}
 -p(1+u^2)^{p-1}   t} {|t|^{1+\alpha}}\,   dt \right| \\
  &\leq & c_5 |x|^{p-\alpha}  \\
  &\leq & c_8 (x^2+a^2)^{(p-\alpha)/2}  .
\end{eqnarray*}
This proves (\ref{eqn:2.13}) for the case of $|x|>a$ and so the
lemma is established. \qed

\bigskip
\noindent
For $\kappa>0$, let $W_t:=\wh S_{\kappa t}$.
It follows from (\ref{eqn:4.6}) that the infinitesimal generator of
$W$ is $\kappa \wh \Delta^{\alpha/2}$ in the sense that for every
$f\in C^2(\R)$,
$$ t\mapsto f(W_t)-\int_0^t \kappa \wh \Delta^{\alpha/2}f(W_s) ds
$$
is a local martingale.

\begin{thm}\label{T:2.8}    Let $f_t(x)$ be the
solution of the backward equation
(\ref{eqn:1.2})
  driven by $W_t:=\wh S_{\kappa t}$. Define $\wt
f_u(z)= f_{\gamma_u} (z)$. Then for every $\alpha, \beta \in (0, 2)$
and $\delta>0$, there is a constant $\kappa =\kappa (\alpha,
\delta)>0$ such that, for every $z=x+iy\in \H$ with $0<y<1$,
 \begin{equation}\label{eqn:2.14}
  \E_z \left[ |\wt f'_u(z)|^\beta ; \, \gamma_u<\infty \right] \leq
e^{-(\beta -\delta) u} \left( x^2 +y^2 \right)^{\beta/2} y^{-\beta}
\qquad \hbox{for } 0<u\leq -\log y.
\end{equation}
\end{thm}

\pf Set
$$ \wt F(u, x, y):=\E_z \left[ |\wt f'_u(z)|^\beta ; \, \gamma_u<\infty \right].
$$
Note that since $\displaystyle \partial_t
f_t(z)=\frac{-2}{f_t(z)-W_t}$, $\displaystyle \partial_t f_t'(z)=
\frac{2f'_t(z)}{(f_t(z)-W_t)^2}$. Thus we have
$$\partial_t \log f_t'(z)=\frac{2}{(f_t(z)-W_t)^2}
$$
and so
$$ \log |f'_t(z)|={\rm Re} \left(\log f_t'(z) \right)=
\int_0^t {\rm Re} \left( \frac{2}{(f_s(z)-W_s)^2}\right) ds
$$
Since $X_t+i Y_t:= f_t (z)-W_t$, it follows that
\begin{equation}\label{deri1}
\log |f'_t(z)|= \int_0^t   \frac{2 {\rm
Re}((X_s-iY_s)^2)}{(X_s^2+Y_s^2)^2} ds =\int_0^t
\frac{2(X^2_s-Y^2_s)}{(X^2_s+Y^2_s)^2}ds
\end{equation}
and so
\begin{equation}\label{deri2}
\log |\wt f'_u(z)|=  \int_0^u \frac{(\wt X^2_s-\wt Y^2_s)}{(\wt
X^2_s+\wt Y^2_s)^2}ds.
\end{equation}
 Thus we have
 \begin{eqnarray}\label{eqn:2.19}
 \wt F(u, x, y) & =& \E_z \left[ \exp \left( \beta \log |\wt f'_u(z)|
\right); \, \gamma_u<\infty  \right] \nonumber \\
 &=&\E_z \left[ \exp \left(
\int_0^u \beta \frac{2(\wt X^2_s-\wt Y^2_s)}{\wt X^2_s+\wt Y^2_s}ds
; \, \gamma_u<\infty\right) \right].
\end{eqnarray}
Observe that by Ito's formula (cf. \cite{HWY}), the infinitesimal
generator $\wt \LL$ of the process $(\wt X, \wt Y)$ is given by
$$ \wt \LL \varphi = -x \frac{\partial \varphi}{\partial x} +y
\frac{\partial \varphi}{\partial y} +   \frac{x^2+y^2}2 \kappa \wh
\Delta^{\alpha/2}_x \varphi,
$$
in the sense that for any $\varphi \in C^2_b(\R^2)$, $\displaystyle
t\mapsto \varphi (\wt X_t, \wt Y_t)-\varphi (\wt X_0, \wt Y_0)
 -\int_0^t \LL \varphi (\wt X_s, \wt Y_s) ds$ is a local martingale.
So formally, when $\alpha \in [1, 2)$, $\wt F$ should satisfy
 \begin{equation}\label{eqn:2.16}
\frac{\partial \wt F}{\partial t}=\wt \LL \wt F+ \beta
\frac{x^2-y^2}{x^2+y^2} \wt F \qquad \hbox{with } \ \wt F(0, x,
y)=1,
\end{equation}
in some sense. Our approach is motivated by this observation.
However (\ref{eqn:2.16}) will not be used in our proof so we can
avoid the delicate questions about the regularity of $\wt F$ and in
which sense the   equation (\ref{eqn:2.16})  holds.

 For $\lambda>0$ and
$\beta>0$, define
$$ \varphi (t, x, y)= e^{-\lambda t} \left( x^2 +y^2
\right)^{\beta/2} y^{-\beta}.
$$
By (\ref{eqn:estimate}), for $\beta>0$, there is a constant
$C_{\beta, \alpha}>0$ such that
$$ | \wh \Delta_x^{\alpha/2} (x^2+y^2)^{\beta/2} |\leq C_{\beta, \alpha}
(x^2+y^2)^{\beta/2 -1} \qquad \hbox{for } |x|\geq 2.
$$
On the other hand, by Lemma \ref{L:2.7},
   there are constants $C_1, C_2>0$, depending only on $\alpha$ and $\beta$,
    such that
$$| \wh \Delta_x^{\alpha/2} (x^2+y^2)^{\beta/2} |\leq C_1
(x^2+y^2)^{(\beta-\alpha)/2} +C_2.
$$
Now take $\beta \in (0, 2)$. We have from above that
$$| \wh \Delta_x^{\alpha/2} (x^2+y^2)^{\beta/2} |\leq c\,
(x^2+y^2)^{\beta/2-1} \qquad \hbox{for } |x| <2 \hbox{ and }
 y\in (0, 1].
$$
 By increasing the value of $C_{\beta, \alpha}>0$ if necessary, we have
\begin{equation}\label{eqn:2.12}
| \wh \Delta_x^{\alpha/2} (x^2+y^2)^{\beta/2} |\leq C_{\beta,
\alpha} (x^2+y^2)^{\beta/2 -1} \qquad \hbox{for every } x\in \R
\hbox{ and } 0<y\leq 1.
\end{equation}
Thus for any $z=(x, y)$ with $x\in \R$ and $0<y\leq 1$, we have
\begin{eqnarray*}
&&  \wt \LL \varphi (t, x, y)+ \beta \frac{x^2-y^2}{x^2+y^2} \varphi (t, x, y) \\
&\leq &\beta  \frac{-x^2+y^2}{x^2+y^2} \varphi (t, x, y)-\beta
\varphi (t, x, y)+  \frac{x^2+y^2}{x^2+y^2} C_{\beta, \alpha} \kappa
\varphi (t, x, y) +\beta \frac{x^2-y^2}{x^2+y^2} \varphi (t, x, y)\\
&=& -(\beta -C_{\beta, \alpha} \kappa )\varphi (t, x, y)
\end{eqnarray*}
So for any $\delta >0$ we can choose $\kappa>0$ small so that
$C_{\beta, \alpha} \kappa <\delta$. Taking $\lambda =\beta -\delta$,
we have
$$ \left(\wt \LL+ \beta \frac{x^2-y^2}{x^2+y^2}\right) \varphi (t, x, y)
\leq -(\beta  -\delta) \varphi (t, x, y) =\frac{\partial}{\partial
t} \varphi (t, x, y)
$$
 for $x\in \R$ and $0<y\leq 1$.
Thus by Ito's formula (cf \cite{HWY}), for each fixed $0<t\leq -\log
y$, with $(\wt X_u, \wt Y_u):=(X_{\gamma_u}, Y_{\gamma_u})$ and
$q(x, y):=\beta \frac{x^2-y^2}{x^2+y^2}$,
$$ M_s:= \varphi (t-s, \wt X_s, \wt Y_s) \exp \left( \int_0^s q(\wt X_u,
\wt Y_u) du\right) {\bf 1}_{\{\gamma_s<\infty\}}
$$
is a supermartingale. It follows that $\E_z M_0\geq \E_z M_t$ and so
$$ \varphi (t, x, y) \geq \E_z \left[ \varphi (0, \wt X_t, \wt Y_t)
\exp \left( \int_0^t q(\wt X_u, \wt Y_u) du\right); \,
\gamma_t<\infty \right].
$$
Since $\varphi (0, x, y) \geq 1 $ for $x\in \R$ and $y\in (0, 1]$
and $\wt Y_t\in (0, 1]$ for every $t\leq -\log y$, we have
$$  \varphi (t, x, y) \geq \E_z \left[
\exp \left( \int_0^t q(\wt X_u, \wt Y_u) du\right); \,
\gamma_t<\infty \right] =\wt F(t, x, y).
$$
This proves the theorem. \qed

\medskip

\section{H\"older continuity}

In this section we will
first
prove that the maps $f_t$ generated by the
truncated stable process $\wh S_{\kappa t}$ are H\"older continuous
a.s.
For small $\kappa$ we obtain explicit estimates for the
exponent.
We will then use Lemma \ref{L:3.1} and the relation between $S$ and $\wh S$
explained below, in order to obtain
H\"older continuity for $f_t$ driven by $S.$

The proof for $\wh S_{\kappa t}$ is similar to the analogous result
for SLE$_\kappa$ with $\kappa\neq4,$ Theorem 5.2 in \cite{RS}. We
begin with an estimate for the derivative $|f_t'|$ of
the backward SLE $\{f_t, t\geq 0\}$ of (\ref{eqn:1.2})
 driven by $W_t=\wh S_{\kappa t}$, using Theorem
\ref{T:2.8}.

\begin{lemma}\label{tderivative}
Let $T>0$.
For $0<\rho<1$ and $\eps>0$
there is $\kappa=\kappa(\rho, \epsilon)>0$ such that
for
$z=x+i y$ with $-R<x<R$ and $0<y<1$, there
is a constant $C>0$ depending on
$T, \alpha, \eps, R$ and $\rho$
so that
 $$
 \P \left( \ \max_{0\leq t\leq T} |f_t'(z)|\geq y^{\rho-1}\ \right) \leq
C y^{2-6\rho -\eps} .
$$
\end{lemma}

\pf Fix
$0\leq t\leq T,$
$z=x+i y$ and write $f_t(z)-W_t = X_t + i Y_t$, $\wt
f_u(z)-W_{\tau_u} = \wt X_u + i \wt Y_u$. Notice $y=Y_0.$
 Recall by (\ref{deri2}),
$$\partial_u \log |\wt f'_u(z)| =
\frac{\wt X_u^2-\wt Y_u^2}{\wt X_u^2+\wt Y_u^2},$$ so that
$$|f_t'(z)| = \exp \left(\int_0^{\log\frac{Y_t}{Y_0}} \frac{\wt X_u^2-\wt Y_u^2}
{\wt X_u^2+\wt Y_u^2} du \right).$$ Let
$$q(u):=\frac{\wt X_u^2-\wt Y_u^2}{\wt X_u^2+\wt Y_u^2}.
$$
Since $|q(u)| \leq 1$, if $Y_t<y^{\rho},$ it follows that
$|f_t'(z)|< y^{\rho-1}.$

\medskip

  On  $\{1\geq Y_t\geq y^\rho \}$, since $q(u) \leq 1$,
 $$
|f_t'(z)| = \exp\left( \int_0^{\log\frac{y^\rho}{y}} q(u) du +
\int_{\log\frac{y^\rho}{y}}^{\log\frac{Y_t}{y}} q(u) du \right)
 \leq |\wt f_{(\rho -1)\log y}' (z)|\,
\frac{Y_t}{y^\rho} \leq |\wt f_{(\rho -1)\log y}' (z)| \, y^{-\rho},
 $$
while on $\{Y_t>1\}$ we have by Lemma \ref{height}
$$|f_t'(z)| = \exp\left(\int_0^{\log\frac{1}{Y_0}} q(u) du +
\int^{\log\frac{Y_t}{Y_0}}_{\log\frac{1}{Y_0}} q(u) du \right) \leq
|\wt f_{-\log y}' (z)| \, Y_t  \leq
(\sqrt{1+4t})
\, |\wt f_{-\log y}' (z)|.
$$
It follows that, on  $\displaystyle \left\{\max_{0\leq t\leq T}
|f_t'(z)| \geq  y^{\rho-1}\right\} \subset \{Y_T\geq y^\rho \}$,
$$\max_{0\leq t\leq T} |f_t'(z)|\leq |\wt f_{(\rho -1)\log y}' (z)| \, y^{-\rho}
1_{\{\gamma_{(\rho-1) \log y} <\infty\}}
+ (\sqrt{1+4t}) \, |\wt f_{-\log y}' (z)| 1_{\{\gamma_{-\log y} <\infty\}}$$
Let $0<\beta<2$ and $\delta>0.$
Then it follows from the above and
 Theorem \ref{T:2.8} that
\begin{eqnarray*}
&& \P\left( \max_{0\leq t\leq T}|f_t'(z)| \geq  y^{\rho-1} \right) \\
 &\leq & \E_z \left[
\left(y^{1-\rho} \max_{0\leq t\leq T} |f'_t(z)|\right)^{\beta};    Y_T \geq y^\rho \right] \\
&\leq & \E_z \left[ (y^{1-\rho})^{\beta} y^{-\rho \beta} |\wt
f_{-\log y}' (z)|^\beta; \gamma_{(\rho-1) \log y} <\infty \right] +
(\sqrt{1+4t})^{\beta}\, \E_z \left[ (y^{1-\rho})^{\beta}  |\wt
f_{-\log y}' (z)|^\beta; \gamma_{-\log y} <\infty \right] \\
&\leq & C  y^{\beta -2\beta \rho} y^{(1-\rho)(\beta -\delta)}
y^{-\beta}+ C y^{\beta -\rho \beta} y^{(\beta-\delta)} y^{-\beta}
\\
&=& C \left( y^{\beta -3\beta \rho +\rho \delta -\delta} + y^{\beta
-\beta \rho -\delta} \right).
\end{eqnarray*}
This establishes the lemma.
 \karo

\bigskip
\noindent
The next theorem says that $z\mapsto f_t(z)$ is
locally uniformly $\gamma$-H\"older
  continuous with any exponent
$\gamma<\frac16$ for $0<\alpha <2$ and small $\kappa.$ We believe
that the correct H\"older exponent is $\frac12$ for every $0<\alpha
<2$. Notice the uniformity in $t\in[0,T],$ which is important later.

\begin{thm}\label{hoelder}
For every $\eps>0$ there is $\kappa>0$ such that with $W_t=\wh
S_{\kappa t}$, for every bounded set $A\subset \H$ and every $T>0,$
a.s. all $f_t$, $0\leq t\leq T$, are H\"older continuous with exponent $1/6-\eps$ on $A$
when $0<\alpha<2$:
$$  |f_t(z)-f_t(z')| \leq C|z-z'|^{\frac16-\eps}
$$
 for all $z,z'\in A$ with a random
constant $C=C(A,\alpha,T,\eps).$
\end{thm}

\noindent {\bf Proof.} Let $R>0$ and $b>0$ be such that $A\subset
[-R,R]\times (0,b].$ It suffices to show that
$$ \max_{0\leq t\leq T}|f'_t(x+i y)| \leq C y^{-\frac56-\eps}
$$
 for all $-R<x<R$ and all
$0<y\leq b.$ By Koebe distortion, it is enough to show this for
dyadic points $z_{j,n} = (j+i) 2^{-n}$, where $n\geq 0$ and $-R 2^n
\leq j \leq R 2^n.$ For every $\eps>0$ there is $\rho> 1/6-\eps$
   such that the exponent $2-4\rho -\eps$ in Lemma \ref{tderivative} is
strictly larger than 1.
  Hence
 \begin{eqnarray*}
\sum_{n=0}^{\infty}\sum_{j=-R 2^n}^{R 2^n} \P\left( \max_{0\leq t\leq T}|f_t'(z_{j,n})|
  >  y^{-\frac56-\eps}\right) <\infty
\end{eqnarray*}
 and the theorem follows as Theorem 5.2 in \cite{RS}. \karo

\medskip
\noindent It immediately follows that for each {\it fixed} $t$, the
map $f_t(z)$ extends continuously to $\overline \H$ a.s. In order to
pass from SLE driven by truncated stable process $\wh S_{\kappa t}$
to SLE driven $S_{\kappa t},$ let's recall the following relation
between $S_t$ and $\wh S_t$. Note that symmetric $\alpha$-stable
process has L\'evy measure $c_\alpha |h|^{-1-\alpha} dh$.
 The jumps
$\{(S_{\kappa t}-S_{\kappa t-}) {\bf 1}_{\{|S_t-S_{t-}|>1\}}, t\geq
0\}$ of size larger than 1 form a Poisson point process with
intensity measure $c_\alpha \kappa  |h|^{-1-\alpha} {\bf
1}_{\{|h|>1\}} dt dh$. The process
$$ \wt S_{\kappa t}:=S_{\kappa t}-\sum_{r\leq t}
     (S_{\kappa r}-S_{\kappa r-}) {\bf 1}_{\{|S_r-S_{r-}|>1\}},
     \qquad t\geq 0,
$$
has the same distribution as $\{\wh S_{\kappa t}, t\geq 0\}$.
 Define $T_0=0$ and
let
$$ T_k:=\inf\{t>T_{k-1}: |S_{\kappa t}-S_{\kappa t-}|>1\} \qquad
\hbox{for } k\geq 1,
$$
be the $k$th jumping time of $S_{\kappa t}$ of size larger than 1.
Then $\{T_k-T_{k-1}, k\geq 1\}$ is a sequence of
i.i.d. exponential
random variables with parameter $\lambda \kappa$. Moreover, the
processes
$$\left\{\{S_{\kappa
(t+T_{k-1})}-S_{\kappa T_{k-1}}, t\in [0, T_k-T_{k-1})\}, k\geq 1
\right\}
$$
 are i.i.d., which are independent copies of $\wh S_{\kappa t}$ killed
at an independent exponential random time $T_1$.
All this tells us that $S_{\kappa t}$ can be constructed as follows.

Let $T_0=0$ and $\{T_k-T_{k-1}, k\geq 0\}$ be an i.i.d. sequence of
exponential random variables with parameter $\lambda \kappa$. Let
$\{\wt S^k_{\kappa t}, t\geq 0\}$ be a sequence of independent
copies of $\{\wh S_{\kappa t}, t\geq 0\}$. Let $\{\xi_k, k\geq 1\}$
be an i.i.d. sequence of random variables with density function
proportional to ${\bf 1}_{|h|>1} |h|^{-1-\alpha}$. These $\{T_k,
k\geq 1\}$, $\{\wt S^k_{\kappa t}, t\geq 1\}$ and $\{\xi_k, k\geq
1\}$ are all independent. For $t>0$, let $n$ be the largest integer
so that $T_n\leq t$. Define
  \begin{equation}\label{eqn:5.1}
    X_t:= \sum_{k=1}^{n-1} \left(\wt S^k_{\kappa (T_k-T_{k-1})}+\xi_k \right)
 +\wt S^n_{\kappa (t-T_n)}.
 \end{equation}
Then $\{X_t, t\geq 0\}$ has the same distribution as $\{S_{\kappa
t}, t\geq 0\}$. From this, we immediately have the following.

\begin{lemma}\label{composition}
For $\kappa >0$, let $\{T_k, k\geq 1\}$, $\{\wt S^k_{\kappa t},
t\geq 1\}$ and $\{\xi_k, k\geq 1\}$ be as in the last paragraph,
which are all independent, and let $X$ be defined by
(\ref{eqn:5.1}). Let $\{f^{(k)}_t, t\geq 0\}$ be SLE driven by
$S^k_{\kappa t}$. For $t>0$, let $n$ be the largest integer so that
$T_n\leq t$. Define
$$ f_t (z):=   \left( f^{(n)}_{t-T_n} (\cdot -X_{T_n})+X_{T_n} \right)\circ
\cdots \circ \left( f^{(2)}_{T_2-T_1}(\cdot -X_{T_1})+X_{T_1}
\right) \circ f^{(1)}_{T_1}(z) .
$$
Then $\{f_t(z), t\geq 0\}$  has the same distribution as the SLE
driven by $W_t=S_{\kappa t}$.
\end{lemma}

\noindent
Because compositions of H\"older continuous maps are H\"older,
from Theorem \ref{hoelder}, Lemma \ref{L:3.1} and Lemma \ref{composition} we obtain the following

\begin{corollary}\label{hoelder2}
For every $0<\a <2$, $\kappa >0$,
 and $W_t = \wh S_{\kappa t}$, for every bounded
set $A\subset \H$ and every $t>0,$ a.s. $f_t$ is H\"older continuous
on $A.$ The same holds for $W_t=S_t$.
\end{corollary}

\section{Hausdorff dimension}

We will now show that the hulls have Hausdorff dimension 1 almost
surely. The situation is similar to  \cite{RS}, Section 8.2: Because
$f_t$ is H\"older continuous, the (box counting) dimension can be
estimated by the convergence exponent of the Whitney decomposition
of $\H\setminus K_t$, which in turn is controlled by the growth of
the derivative $f_t'$ towards the boundary $\R$ of $\H.$
For a Borel set $K\subset \R^2$, we use $\dim_H K$ to denote its
Hausdorff dimension.

\begin{thm}\label{dimension}

\noindent For each $0<\alpha<2$, $\kappa >0$,  and $W_t = S_{\kappa
t}$ (or $W_t = \wh S_{\kappa t}$)
$$\dim_H K_t = 1$$
for all $t\geq 0$, almost surely.

\end{thm}

\noindent
Since $K_t$ has empty interior by \cite{GW}, $K_t\cup\R = \partial
(\H\setminus K_t).$ Because $g_t^{-1}$ has the same distribution as
$f_t$ (for fixed t), it thus suffices to show that the boundary of
$H_t = f_t(\H)$ has dimension 1 a.s. Denote by $N(\eps)=N(\eps,A)$
the minimal number of disks of radius $\eps$ needed to cover a set
$A\subset\C$.

The following is an analog of the upper bound for the dimension of the outer SLE boundary,
Theorem 8.6 in \cite{RS}.

\begin{thm}\label{cover}
\noindent For each $0<\alpha<2$ and $1<a<2$, there is $\kappa>0$
such that with $W_t=\wh S_{\kappa t}$, for all $T>0, h>0$ and $R>0$,
a.s. we have

$$\lim_{\eps\to0}
\eps^a\,\max_{0\leq t\leq T}\, N\bigl(\eps, f_t[-R, R] \cap \{y>h\}\bigr)= 0.
$$

\end{thm}

\noindent {\bf Proof.} As in \cite{RS}, Section 8.2, consider a
Whitney decomposition of $H_t$ (that is a covering of $H_t$ by
essentially disjoint closed squares $Q\subset H_t$ with sides
parallel to the coordinate axes such that the side length $d(Q)$ is
comparable  to the distance of $Q$ from the boundary of $H_t$, and
such that $d(Q)$ is an integer power of 2). Denote by $W_t$ the
collection of those squares $Q$ for which $Q \cap f_t([-R, R]\times
(0,\infty))\cap \{y>h\} \neq\emptyset$, and let
$$S(a)= \max_{0\leq t\leq T} \sum_{Q\in W_t} d(Q)^a \leq \infty.$$
Then the proof of Theorem 8.6 in \cite{RS} (the last displayed formula) shows
that, for each $0\leq t\leq T,$
$$
N\bigl(2^{-n}, f_t[-R, R]\cap \{y>2\,h\}\bigr)
\le C(\omega)\, 2^{(n+O(\log n))\,a} S(a)\,.
$$
The factor $C(\omega)$ comes from the H\"older norm of $f_t$ and is random,
but does not depend on $t$ or $n.$
The theorem follows at once if we show $S(a)<\infty$ a.s. To this end, we will show
$$\E[S(a)] < \infty$$ for $a>1,$ in analogy with the upper bound in Theorem 8.3 of \cite{RS}.
By the Koebe distortion theorem, again writing $z_{j,n} = (j+i) 2^{-n}$, the quantity
$$\wt S(a) = \max_{0\leq t\leq T}\,\sum_{n=0}^{\infty}\sum_{j=-R 2^n}^{R 2^n}
\1_ {\{\Im f_t(z_{j,n})>h\}}
\, |f_t'(z_{j,n})\, 2^{-n}|^a $$
is comparable to $S(a)$ (see (8.2) and Lemma 8.4 of \cite{RS} for
the details), in particular $S(a)\leq C \wt S(a)$ for some universal
$C.$
For $0\leq t\leq T,$ Lemma \ref{height} yields
$$\1_ {\{\Im f_t(z_{j,n})>h\}} |f_t'(z_{j,n})|
\leq \1_{\{\gamma_{\log{(h/y_{j,n})}}  <\infty\}} |f_t'(z_{j,n})|
\leq C_{T,h}
|\wt f_{\log (2^n h)}' (z_{j,n})| {\bf 1}_{\{\gamma_{\log{(2^nh)}}
<\infty\}} .$$
Now Theorem \ref{T:2.8} with
$\beta=a$ and $\delta< a-1$ implies
$$\E[S(a)] \leq \sum_{n=0}^{\infty} 2 R 2^n 2^{-a n} 2^{n \delta} < \infty,$$
thus proving the theorem.
\karo

\bigskip\noindent
The next lemma says that the Hausdorff dimension of boundaries of simply connected domains
does not increase under finite composition.

\begin{lemma}\label{dim-composition}
Let $1<a<2,$ let $f^{(j)}:\H\to\H\setminus K^{(j)},$ $1\leq j\leq
n,$ be conformal maps, and let $f=f^{(n)}\circ \cdots \circ
f^{(1)}.$ If $\dim_H K^{(j)}\leq a$ for all $j,$ then $\dim_H
\partial f(\H)\leq a.$
\end{lemma}
\noindent

\pf
For the proof, just notice that
$$\partial f(\H) =
K^{(n)}\cup f^{(n)}(K^{(n-1)}) \cup \cdots \cup f^{(n)} (f^{(n-1)}
(\cdots (f^{(2)}(K^{(1)}))))$$
and that each of the sets in the union has dimension $\leq a$ because each map $f^{(j)}$
is smooth in $\H.$\qed

\bigskip
\noindent {\bf Proof of Theorem \ref{dimension}.} Let $1<a<2$ and
let $\kappa>0$ be as in Theorem \ref{cover}. Let $f_t$ be driven by
$S_{\kappa t},$ and factor $f_t$ according to Lemma
\ref{composition}. Then by Theorem \ref{cover}, the hulls of the
factors of $f_t$ have Hausdorff dimension $\leq a,$ and thus $\dim_H
\partial f_t(\H)\leq a$ by Lemma \ref{dim-composition}. Letting $a$
tend to 1, we see that the hulls driven by $W_t=S_{\kappa t}$ have
Hausdorff dimensional at most 1.
Because the boundary of the simply connected domain $\H\setminus K_t$
is connected, and because $K_t\cap \H\neq\emptyset$, we have
$\dim_H K_t\geq 1$ and conclude $\dim_H K_t= 1$ for every $t>0.$
By the scaling Lemma \ref{L:3.1}, the hulls driven by $W_t=S_t$ have
the same dimension as the hulls of $S_{\kappa t}$. Finally, the
hulls of $\wh S_{\kappa t}$ for an arbitrary (not neccessarily
small) $\kappa$ can be recovered from the hulls of $S_{\kappa t}$,
and therefore have dimension 1, by removing the jumps of $W_{\kappa
t}$, similar to Lemma \ref{composition}. \qed

\section{Trace continuity}\label{cont}

The purpose of this section is to prove the following
\begin{thm}\label{continuity}
Fix $\alpha\in(0,2)$ and $\kappa >0$. Let $W_t=S_{\kappa t}$ or $W_t= \wh S_{\kappa t}$.
 Then almost surely,
 for each $t>0$ the limit
$$\gamma(t) = \lim_{z\to W_t; z\in\H} g_t^{-1}(z)$$
exists, the function $t\mapsto \gamma(t)$ is RCLL, and
$K_t = \overline{ \gamma[0,t]}$.
\end{thm}

\bigskip \noindent
From Theorem \ref{hoelder} we know that for each fixed $T,$ $f_T(z)$
extends continuously to $\overline \H$ a.s. Because the hulls $K_T$
have the same law as $\H\setminus f_T(\H)-W_T,$ they are locally
connected a.s. In general, this does not imply that the subsets
$K_t\subset K_T$ for $t<T$ are locally connected too
(for instance,
it is possible that $K_t$ is not locally connected at some time $t_0$,
but that due to "swallowing" $K_{t_0}$ is contained in the interior
of $K_t$ for some $t_1>t_0,$ and that the boundary of $K_{t_1}$ is
smooth).  Nor does the
equicontinuity of $f_t(z)$ generally imply equicontinuity of
$g_t^{-1}(z).$ For instance, if $W_t=c\sqrt{1-t}$ with $c=2\sqrt{3},$
then $f_t(z)$ is equicontinuous ($\H\setminus f_t(\H)$ is a halfdisc
of radius proportional to $\sqrt{t}$) whereas $g_t^{-1}$ is not
($K_t$ is an arc of a semicircle up to time 1 when $K_t$ is a
semidisc).
Because of the tree structure of the hulls, our
situation is better:

\begin{prop}\label{equicont}
Let $W_t=S_{\kappa t}$ or $W_t= \wh S_{\kappa t}$. For each
$0<\alpha<2$ and each $T>0,$ a.s. each of the maps $g_t^{-1}$,
$0\leq t\leq T$,  has a continuous extension to $\overline\H$ (which
we again denote $g_t^{-1}$). Moreover, the maps $\{g_t^{-1}$, $0\leq
t\leq T \}$, are equicontinuous on $\overline\H$.
\end{prop}

\bigskip
\noindent
We postpone the proof until the end of this section and continue with the
\bigskip

\noindent {\bf Proof of Theorem \ref{continuity}.}
Fix $\alpha\in(0,2)$ and $T>0,$ and let $t\leq T.$ Because $g_t^{-1}$ has a
continuous extension to $\overline \H$ by Proposition \ref{equicont},
$\gamma(t) = \lim_{z\to W_{t}; z\in\H} g_t^{-1}(z)$ exists,
and $\gamma (t)=g_t^{-1}(W_t).$
The equicontinuity of $g_t^{-1},$ together with the pointwise continuity of
$t\mapsto g_t^{-1}(z)$, easily implies the continuity of $(t,z)\mapsto g_t^{-1}(z)$
on $[0,T]\times\overline\H.$ It follows immediately that $\gamma$ is RCLL.

\noindent
To prove $K_t = \overline{\gamma[0,t]}$, first let $z=\gamma(t)=g_t^{-1}(W_t).$
Clearly $T_z\leq t$, so that $z\in K_t$ (see Section \ref{S:3.1} for the notation).
Because $K_t$ is closed, we have $K_t \supset \overline{\gamma[0,t]}.$
Conversely, if $w\in K_t,$ then $T_{w}\leq t.$ Then either
$\displaystyle \liminf_{s\to T_{w}-} |g_s(w)-W_s| = 0,$ and the
continuity of $(s,z)\mapsto g_s^{-1}(z)$ implies $\displaystyle
\liminf_{s\to T_{w}-} |g_s^{-1}(g_s(w))-g_s^{-1}(W_s)| = 0,$ which
yields $w\in \overline{\gamma[0,t]}.$ Or $\displaystyle
\liminf_{s\to T_{w}-} |g_s(w)-W_s| > 0$ and
$g_{T_{w}}(w)=W_{T_{w}}$, which means $w=\gamma(T_{w})\in
\overline{\gamma[0,t]}$. It follows that $K_t \subset
\overline{\gamma[0,t]}.$ \qed

\bigskip
\noindent In order to prove Proposition \ref{equicont}, we need a
variant of a theorem of Warschawski \cite{W} about the modulus of
continuity of conformal maps of the disc. Roughly speaking, after
suitable normalization the modulus only depends on the "roughness"
of the boundary of the domain as measured by the size of
bottlenecks. Let $G\subset\wh \C$ be a simply connected domain and
$a\in G$ be a marked point (in \cite{W}, $a=0$ whereas here we will
have $p=\infty$). A {\it crosscut} of $G$ is a simple arc
$\{\sigma(t),$ $0\leq t\leq 1\}$ that lies in $G$ except for the
endpoints $\sigma(0),\sigma(1)\in\partial G.$ Every crosscut
separates $G$ into two connected components. If $a\notin \sigma,$
denote $G(\sigma)$ the component that does not contain $p$ in its
closure. Following Warschawski, define
$$\eta_G(\delta) = \sup_{\diam \sigma \leq \delta} \diam G(\sigma).$$
Thus $\eta_G(\delta)
 \to 0$ as $\delta\to0$ if and only if $\partial G$
is locally connected. Now assume that $G=\H\setminus K$ and that
$f:\H\to G$ is the hydrodynamically normalized conformal map, $f(z)
= z + a/z + O(1/z^2)$ near $\infty.$  Denote
$$
 \omega_f(r) = \sup\left\{|f(z)-f(z')|:z,z'\in\H \hbox{ with }  |z-z'|\leq r \right\}
$$
 the modulus of
continuity of $f.$ The following is Theorem I of Warschawski
\cite{W}, except for the different normalization. His proof carries
over with only minor modifications.

\bigskip
\noindent
\begin{thm}\label{warschawski}
For each $R>0$ and each function
$\eta(\delta)$ with $\eta(0+)=0$
there is a function
$\omega(r)$ with $\omega(0+)=0$ such that the following holds:
If $K\subset \{|z|<R\},$  and if
$\eta_G(\delta)\leq \eta(\delta)$ for all $\delta,$ then
$$\omega_f(r) \leq \omega(r)$$
for all $r>0.$
\end{thm}

\bigskip
\noindent {\bf Proof of Proposition \ref{equicont}}. Fix $T>0.$
Because $f_T(z)$ extends continuously to $\overline \H$ a.s. by
Corollary \ref{hoelder2}, and because the hulls $K_T$ have the same
law as $\H\setminus f_T(\H)-W_T,$ we have
$$\eta_{G_T}(0+)=0$$
a.s., where $G_T=\H\setminus K_T.$ By Theorem 1.3 (i) in \cite{GW},
we know that $K_T$ and hence $K_t,$ $0\leq t\leq T,$
does not have interior points. Hence every crosscut $\sigma$ of
$\H\setminus K_t$
  can be decomposed into crosscuts $\sigma_j$ of
$\H\setminus K_T$ such that $(\H\setminus K_t)(\sigma) \subset
\overline{\cup_j  (\H\setminus K_T)(\sigma_j)}.$
It follows that
$\eta_{G_t}(\delta)\leq \delta + 2\eta_{G_T}(\delta),$ for all $t\leq T.$ Now Proposition
\ref{equicont} follows from Theorem \ref{warschawski}. \karo

\end{doublespace}

\end{document}